\documentclass[11pt]{article}
\usepackage{graphics, amsfonts, amsmath, amsthm, amstext, amssymb, amscd, color, rotating, makeidx}
\usepackage{calligra,mathrsfs} 
\usepackage{ulem} 
\usepackage{tikz-cd} 
\usepackage[square,sort,comma,numbers]{natbib}
\usepackage{hyperref}
\hypersetup{
    colorlinks, 
    citecolor=red,
    filecolor=black,
    linktoc=all,     
    linkcolor=blue, 
    urlcolor=blue 
}
\usepackage[margin=0.9in]{geometry} 
\usepackage{fancyhdr}  
\newtheorem{theorem}{Theorem}[section]

\newtheorem{lemma}[theorem]{Lemma}
\newtheorem{proposition}[theorem]{Proposition}

\newtheorem*{remark}{Remark}

\newcommand{\Aut}{\operatorname{Aut}}
\newcommand{\barD}{\overline{\partial}}
\newcommand{\ddbar}{\partial\overline{\partial}} 
\newcommand{\Vol}{\operatorname{Vol}} 
\newcommand{\Hol}{\operatorname{Hol}} 
\newcommand{\mult}{\operatorname{mult}} 
\newcommand{\sech}{\operatorname{sech}} 
\newenvironment{claim}[1][]{\par\smallbreak
   \noindent \underline{\textsc{Claim}}. \rmfamily}{\smallbreak} 
\newenvironment{cproof}[1][Proof]{%
  \proof[#1]%
  }{\endproof}

\usepackage[]{amsmath, amsthm, amsfonts, verbatim, amssymb,pb-diagram, stmaryrd, wasysym,tikz}
\usepackage{epsfig}

\def\oM{\overline M}

\def\bC{{\mathbb C}}
\def\bN{{\mathbb N}}

\def\tM{\widetilde M}

\def\ms{\medskip}

\def\pd{\partial}

\def\oz1{d{\overline z}^1}
\def\oz2{d{\overline z}^2}
\def\oz3{d{\overline z}^3}

\def\oI{\overline I}

\def\oz{\overline z}

\def\op{\overline D}
\def\oIq1{\oI_1\cdots\oI_{q-1}}
\def\oIq2{\oI_1\cdots\oI_{q-2}}
\def\op{\overline \partial}

\def\Aut{\mbox{Aut}}

\def\ddz{\frac{\pd}{\pd z}}

\def\ddw{\frac{\pd}{\pd w}}
\def\ddz1{\frac{\pd}{\pd z^1}}
\def\ddw1{\frac{\pd}{\pd w^1}}
\def\oddw1{\overline{\frac{\pd}{\pd w_1}}}

\def\tomega{\widetilde{\omega}}

\def\tw{\widetilde w}
\def\tu{\widetilde u}
\def\tR{\widetilde R}

\def\mult{\mbox{mult}}

\begin{document}

\title{Carath\'eodory hyperbolicity, volume estimates and level structures over function fields} 

\author{
Kwok-Kin Wong\footnote{
School of Mathematical Sciences, Shenzhen University, Guangdong 518060, People’s Republic of China.
Email: kkwong@szu.edu.cn} \quad
and 
\quad 
Sai-Kee Yeung\footnote{
Department of Mathematics, Purdue University, 150 N. University Street, West Lafayette, IN 47907-1395, USA.
Email: yeungs@purdue.edu}
}
\date{}
\maketitle

\begin{abstract}  
{\it  
We give a generalization of the nonexistence of level structures as in \cite{Nade1989}, \cite{Nogu1991} and \cite{HT2006} for quasi-projective manifolds uniformized by strongly Carath\'eodory hyperbolic complex manifolds. Examples include moduli space of compact Riemann surfaces with a finite number punctures and locally Hermitian symmetric spaces of finite volume.  This leads to the nonexistence
of a holomorphic map from a Riemann surface of fixed genus into the  compactification of such a quasi-projective manifold when the level structure is sufficiently high.
To achieve our goal, we have also established some volume estimates for mapping of curves into these manifolds, extending some earlier result of \cite{HT2006}
to a more general setting.  A version of Schwarz Lemma applicable to manifolds equipped with nonsmooth complex Finsler metric is also given. 
}
\end{abstract}

\bigskip

\noindent \textbf{Acknowledgments}. The authors would like to thank the referees for valuable remarks and suggestions which lead to improvements of the exposition.

\section*{Introduction}
A complex manifold 
 is said to be
strongly Carath\'eodory hyperbolic if and only if its infinitesimal Carath\'eodory pseudo-metric $g_C$ is nondegenerate and its Carath\'eodory pseudo-distance function $d_{C,M}$ is complete nondegenerate 
(for detail definitions, c.f. \cite{WY2022} or \cite{Koba1998}). 
These manifolds are known to have many nice hyperbolicity properties. As a continuation of \cite{WY2022}, we study in \cite{WY2023} complex hyperbolicity on a tower of
quasi-projective manifolds, which include in particular the cases of genus $0$ and $1$ algebraic curves.  In this paper, we study algebraic curves of genus at least $2$ in quasi-projective manifolds that are uniformized by strongly Carath\'eodory hyperbolic complex manifolds, which requires a completely different set of techniques comparing to genus $\leqslant 1$ cases
as in \cite{WY2023}.
The arguments in this paper make use of the results of \cite{WY2022} and is parallel to those in \cite{WY2023}. 

A manifold $M$ is said to support a tower of coverings $\{M_i\}_{i=1}^\infty$ with $M_1=M$  if for each $i$, there is a finite unramified covering $M_{i+1}\rightarrow M_{i} $ such that the fundamental groups  $\pi_1(M_{i+1})\triangleleft \pi_1(M_1)$ as normal subgroup of finite index, and that $\cap_{i=1}^\infty \pi_1(M_{i})=\{1\}$.

Our first main result is a generalization of \cite{Nade1989}, \cite{Nogu1991} and \cite{HT2006}:
\begin{theorem}
Let $M=\oM- D$ be a quasi-projective manifold uniformized by a strongly Carath\'eodory hyperbolic manifold $\tM$, or a Carath\'eodory hyperbolic manifold equipped with a smooth bounded
plurisubharmonic exhaustion function. 
Suppose $M$ supports a tower of coverings  $\{M_i\}_{i=1}^\infty$. Assume the following properties hold:
\vspace{-0.2cm}
\begin{enumerate}
\item[(i).] There exists a complete K\"ahler-Einstein metric $g_{KE}$ negative scalar curvature on $\tM$ satisfying
\begin{equation}
 c\cdot g_C\leqslant g_{KE}\leqslant \frac1{c}\cdot g_C
\label{CKEest} 
\end{equation}
for some constant $c>0$.
\item[(ii).] 
For $i\geqslant 0$, $M_{i+1}\rightarrow M_{i}$ extends to a finite ramified covering $\oM_{i+1}\rightarrow \oM_{i}$ between projective manifolds, where the boundary divisor $D_i=\oM_i- M_i$ is of simple normal crossings.
\end{enumerate}
\vspace{-0.3cm}
Let  $g_0\geqslant 2$
be a fixed nonnegative integer.  Then there exists $i_0\geqslant 0$ such that for $i\geqslant i_0$, and any Riemann surface $S$ of genus $g_0$, any non-constant holomorphic map $f:S\rightarrow \oM_i$ has image $f(S)\subset D_i$.
\label{MT1}
\end{theorem}
\vspace{-0.4cm}

\begin{remark}
(a) In condition (i), if the projective compactifications $\overline{M_i}$ of $M_i$'s are not smooth, then one may take a resolution of singularities $M_i'$ of $M_i$ and apply the result for the smooth case to see that the image $f(S)$ lifts to lie in $M_i'-M_i$, which blows down to $D_i=\overline{M_i}-M_i$. Hence the conclusion also follows.

\ms
\noindent
(b) If the fundamental group $\pi_1(M)$ is residually finite, i.e. the intersection of all normal subgroups of finite index of $\pi_1(M)$ is trivial, then we know that $M$ supports a tower of covering. In general we don't know if the strong Carath\'eodory hyperbolicity would imply the residually finiteness of the fundamental group $\pi_1(M)$.

\ms
\noindent
(c)
As is considered in \cite{Deng2023}, if $\pi_1(M)$ admits a linear representation $\rho$ to a reductive group, then there is an induced mapping from $M$ to a locally symmetric space. If the image $\rho(\pi_1(M))$ is non-trivial or finite, then $\rho(\pi_1(M))$ will automatically be residually finite and hence corresponding to a tower.
\end{remark}

Note that the existence of such $f$ corresponds to existence of a level structure over the function field of $S$ as a projective algebraic curve.  The statement of the Theorem
can be understood as a statement that $M$ has no level structure over a function field as above when the level is sufficiently high.

Examples satisfying the conditions in Theorem \ref{MT1} include finite volume quotients of HHR/uniform squeezing domains introduced in \cite{LSY2004}  and \cite{Yeun2009}, see also \cite{Yeun2005}, which on the other hand include examples 
of a moduli space of compact Riemann surface with a finite number of punctures and any locally Hermitian symmetric spaces $M_{LS}$.  The fact that the Carath\'eodory distance is complete for HHR/uniform squeezing domains 
can be found in \cite{Yeun2022}.

\begin{remark}
Suppose $M$ is a bounded domain in $\mathbb{C}^n$ with complete Carath\'eodory distance, then $M$ is complete K\"ahler hyperbolic \cite[Theorem 1]{Yeun2022}. The same proof applies in
case that $M$
 is a strongly Carath\'eodory hyperbolic complex manifold.
\end{remark}

Our second main result is confined to Hermitian locally symmetric space or a moduli space of Riemann surfaces: 

\begin{theorem} 
 Let  $g_0\geqslant 2$.\\
1). Let $M=M_{LS}$ be a Hermitian locally symmetric manifold of finite volume and complex dimension at least $2$. Write $M=M_{0}$. 
\vspace{-0.2cm}
\begin{enumerate}
\item[(a).] If $M_{LS}$ is arithmetic, consider $\oM_{LS}$ the Baily-Borel-Satake compactification of $M_{LS}$, and a tower of coverings $\{M_i\}_{i=0}^\infty$ coming from a level structure
; or 
\vspace{-0.3cm}
\item[(b).] If $M_{LS}$ is a nonarithmetic complex ball quotient, consider $\oM_{LS}$ the Siu-Yau compactification of $M_{LS}$ and a tower of coverings $\{M_i\}_{i=0}^\infty$.
\end{enumerate}
\vspace{-0.3cm}
Then there exists $k_o\geqslant 0$ sufficiently large such that if $k\geqslant k_o$, $\overline{M_k}$ does not contain an algebraic curve of genus $\leqslant g_0$.

\noindent 2). Let $M_{g}$  be the moduli space of Riemann surfaces of genus $g\geqslant 2$ and  $\overline{M_g}$ be the Deligne-Mumford compactification of $M_g$.  Write $M^0_g=M_g$. Let
 $\{M^k_g\}_{k=0}^\infty$ be a tower of coverings. Then there exists $k_o\geqslant 0$ sufficiently large such that if $k\geqslant k_o$, we have:
 \vspace{-0.2cm}

\begin{enumerate}
\item[(I).] $\oM_g^k$ does not contain any embedded curves of genus $g_0$ intersecting non-trivially with $M_g^k$; and
\vspace{-0.3cm}
\item[(II).] The set of genus $g_0$ curves on $\pd\oM^k_g$ are contained in 
 the fiber of  a puncture forgetting projection map on an irreducible component of a stratum of the boundary $\pd \oM^k_g$.
  
  \end{enumerate}
\vspace{-0.3cm}  
 \label{MT2}
\end{theorem}

Similar to the earlier remarks after Theorem \ref{MT1},  related to \cite{Nade1989}, \cite{Nogu1991} and \cite{HT2006},
the statement of Theorem \ref{MT2}
can be understood as a statement that the compactification $\oM$ of $M$ has no level structure over a function field when the level is sufficiently high
in the cases of Baily-Borel-Satake or Siu-Yau compactification of a locally Hermitian symmetric spaces. For the Deligne-Mumford compactification of moduli of Riemann surfaces of genus $g\geqslant 2$, the possible appearance of a curve of genus $g_0$
can be read off from the discussions in the proof of Theorem \ref{MT2}. See also \cite{Bru} for some related results, for which the authors thank Soheil Memariansorkhabi for bringing to their attention
after the acceptance of the paper.

\bigskip

The proof of Theorem \ref{MT1} depends on a crucial estimate in the volume of a curve in the coverings of the manifold.  Due to its relation to other
problems of interests, we state it as separate result.
\begin{theorem}
Let $M=\oM- D$ be a quasi-projective manifold supporting a tower of coverings as given in Theorem \ref{MT1}.  Let $S$ be a compact Riemann surface genus $g_S\geqslant 2$ and $f:S\rightarrow \overline{M}_i$ be a nonconstant holomorphic map. Write $V=f(S)$. Suppose $V\cap M_i\neq \emptyset$.   
 Let $x\in V\cap M_i$, we consider the volume $\Vol_{_{KE}}(V\cap B_{_{KE}}(x;R))$ with respect to $g_{_{KE}}$, where $B_{_{KE}}(x;R)\subset M_i$ is the geodesic ball of radius $R>0$ centred at $x$ with respect to $g_{_{KE}}$. Write $\pi:\tM\rightarrow M_i$ as the universal covering map. Let $\widetilde{x}\in \tM$ be such that $\pi(\widetilde{x})=x$.
 $\widetilde{V}\subset \tM$ is a local lifting of $V\cap M_i$ near $x$Then given any $L>0$. For $R$ sufficiently large, $\Vol_{_{KE}}(\widetilde{V}\cap B_{_{KE}}(\widetilde{x};R))\geqslant L$ on $\widetilde{M}$.
\label{MT3}
\end{theorem}

For Hermitian locally  symmetric spaces, it is observed in  \cite{HT2006}
that Theorem \ref{MT3} is crucial for a proof of Theorem \ref{MT1}.   The proof of Theorem \ref{MT3} for locally
Hermitian symmetric spaces was given in  \cite{HT2006}.  The proof of Theorem \ref{MT3} for horizontal slices of period domains 
 was given in \cite{BT2019}.  Both of these results play significant roles in the functional transcendence, such as
 the Ax-Lindemann or, more generally, Ax-Schanuel type problems, see for example \cite{KUY2016, BT2019,MPT2019}.  The results of \cite{HT2006} makes use of
 the non-positivity in the Riemannian sectional curvature of the Bergman metric in the Hermitian symmetric spaces, and the results of \cite{BT2019}
makes use of the
 fact that on the horizontal slices of the period domains, the natural Hermitian metric has holomorphic sectional curvature bounded from above by a negative constant.
 Our contribution is to show that appropriate volume estimates still hold when such strong metric 
 and curvature conditions were not available, replacing by milder conditions from the perspective of Carath\'eodory metric.

 For the proof of Theorem \ref{MT1}, we have also established a version of Schwarz Lemma (Lemma \ref{schwarz})  in which the target is a not necessarily smooth complex Finsler manifold, a result that we need but cannot find in the literature.
The other ingredients include a general deformation theoretic argument taken from \cite{HT2006},
which guarantees an injectivity radius lower bound of curves (Proposition \ref{injlowbdd}),
and the argument to exploit the ramification at boundary divisors along the tower of coverings \cite{Nade1989,Nogu1991}. The essential observation here is that these methods are general enough to be applicable to our case in view of our previous work \cite{WY2022}.

 Theorem \ref{MT2} illustrates possible applications of Theorem \ref{MT1} by two series of examples, namely the Hermitian locally symmetric spaces and the moduli space of Riemann surfaces.  The proof of  Theorem \ref{MT2}
follows from an appropriate understanding of the structure of strata in the compactification and repeated
iteration of Theorem \ref{MT1}.

The organization of the article is as follows. In \S\ref{crs}, some basic facts of compact Riemann surface are recalled in order to fix notations. In \S\ref{vol}, the volume estimates are derived, first with the lower estimates and then the upper estimates. The argument for lower bound of injectivity radius and ramification are given respectively in \S\ref{inj} and \S\ref{ramind}. The proof of Main Theorem \ref{MT1} and \ref{MT2} are given in \S\ref{pfMT}.

\bigskip

\section{Compact Riemann surface}
\label{crs}
Let $S$ be a compact Riemann surface of genus $g$. Suppose $h$ is a complete Hermitian pseudo-metric on $S$. In local coordinate $z$, $h=2\lambda dz\otimes d\overline{z}$ with the associated K\"ahler form $\omega_h=\lambda \frac{\sqrt{-1}}{2\pi}dz\wedge d\bar{z}$. The Gaussian curvature of $h$ is given by
\[
K_h:= -\frac{1}{4}\frac{\Delta \log \lambda}{\lambda},
\]
where the Laplacian $\Delta =4\frac{\partial^2}{\partial z\partial\bar{z}}$ is interpreted in the sense of distribution in general. Let $\chi(S)$ be the Euler characteristic of $S$. If $S$ is of genus $g=g(S)$, then $\chi(S)=2-2g$. By Gauss-Bonnet Theorem,
\[
\int_S K_h\omega_h = \chi(S)=2-2g.
\]
On the other hand, 
the holomorphic tangent bundle $TS$ is a line bundle equipped with metric $h$. The first Chern class of $S$ is given by 
\[
c_1(S)=c_1(TS)=-\frac{\sqrt{-1}}{2\pi}\ddbar \log \lambda=K_h \omega_h.
\]
It follows that the degree of line bundle
\[
\deg (TS) := \int_R c_1(TS) = 2-2g.
\]

\section{Volume estimates}
\label{vol}
In the following, we adopt the notations as in the assumptions of Theorem \ref{MT1} and \ref{MT2}. 

On the universal covering $\tM$ of $M$, there is a K\"ahler-Einstein metric $g_{_{KE}}$ of negative scalar curvature, which is unique if the Einstein constant is fixed.
From uniqueness of the K\"ahler metric with fixed Einstein constant as above, $g_{_{KE}}$ is invariant under $\Aut(\tM)$ so that it descends to $M$. 
Since both $g_C$ and $g_{_{KE}}$ are invariant under $\Aut(\tM)$, they descend to $M_i$ and will be denoted by the same notations.

\subsection{Lower bound in volume}
Let $S$ be a compact Riemann surface genus $g\geqslant 2$ and $f:S\rightarrow \overline{M}_i$ be a nonconstant holomorphic map. Write $V=f(S)$. Suppose $V\cap M_i\neq \emptyset$.   
 Let $x\in V\cap M_i$, we consider the volume $\Vol_{_{KE}}(V\cap B_{_{KE}}(x;R))$ with respect to $g_{_{KE}}$, where $B_{_{KE}}(x;R)\subset M_i$ is the geodesic ball of radius $R>0$ centred at $x$ with respect to $g_{_{KE}}$. Write $\pi:\tM\rightarrow M_i$ as the universal covering map. Let $\widetilde{x}\in \tM$ be such that $\pi(\widetilde{x})=x$.
Denote by $\rho_i(x)$ the injectivity radius at $x\in M_i$ with respect to $g_{_{KE}}$. 
For $0<R<\rho_i(x)$, there is a biholomorphic isometry $\pi: (B_{_{KE}}(\widetilde{x},R), g_{_{KE}})\rightarrow (B_{_{KE}}(x,R),g_{_{KE}})$. To find $\Vol_{_{KE}}(V\cap B_{_{KE}}(x;R))$, it suffices to find $\Vol_{_{KE}}(\widetilde{V}\cap B_{_{KE}}(\widetilde{x};R))$, where $\widetilde{V}\subset \tM$ is a local lifting of $V\cap M_i$ near $x$. For simplicity, from now on we write $\widetilde{V}=V$ and $\widetilde{x}=x$.

\begin{proposition}
Given any $L>0$. For $R$ sufficiently large, $\Vol_{_{KE}}(V\cap B_{_{KE}}(x;R))\geqslant L$ on $\widetilde{M}$.
\label{lvol}
\end{proposition}

\noindent For the proof of Proposition \ref{lvol}, we start with some preparations.

Let $R>0$. Denote by $[V]$ the closed $(n-1,n-1)$-current corresponding to $V$. For the purpose of obtaining a lower estimate of the volume 
\[
\Vol_{_{KE}}(V\cap B_{KE}(x;R)) = \int_{B_{_{KE}}(x;R)} [V]\wedge\omega_{_{KE}},
\]
it suffices to consider a smooth point $x\in V$.
From assumption (i) in Theorem \ref{MT1},  $B_{g_C}(x;c R)\subset B_{{_{KE}}}(x;R)\subset B_{g_C}(x;\frac{1}{c} R)$. The distance function obtained by integrating $g_C$ is exactly the inner distance function with respect to $d_C$ and thus must always $\geq d_C$ (c.f. \cite[Theorem 4.2.7]{Koba1998}). Therefore 
\begin{equation}
B_{d_C}(x;cR)\subset B_{g_C}(x;cR).
\label{Bsubset}
\end{equation}

For $z\in \Delta$, the infinitesimal Poincar\'e metric is defined by $ds^2_\Delta = \frac{dz\otimes d\overline{z}}{1-|z|^2}$. The corresponding Poincar\'e
distance between $0$ and $z$ is given by $\ell_P(z):=\ell_P(0,z)=\frac12\log\frac{1+|z|}{1-|z|}$. Consider the Carath\'eodory distance 
\[
d_C(x,y)=\sup \{ \ell_P(f(x),f(y)) \mid f\in \Hol(M,\Delta)\} ,\quad x,y\in \tM.
\]
Let $x\in\tM$ be fixed. By homogeneity of $\Delta$ and Arzela-Ascoli Theorem, $d_C(x,y)=\ell_P(h(y))$ for some holomorphic map $h:\tM\rightarrow \Delta$ such that $h(x)=0$. For $y\in \tM$, write 
\[
\ell_C(y):=d_C(x,y).
\]
We define
\begin{equation}
r_C(y)=\tanh(\ell_C(y)), \quad y\in \tM.
\end{equation}
Note that for fixed $x\in \tM$, 
\[
r_C(y)
=\tanh\bigg(\frac{1}{2}\log \frac{1+|F(y)|}{1-|F(y)|}\bigg)
=|F(y)|, \quad y\in \tM,
\] 
for some $F\in \Hol(\tM,\Delta)$ such that $F(x)=0$. It follows  that 
\[
r_C(y)= \sup \{ |f(y)| : f\in \Hol(M,\Delta), f(x)=0\}.
\]
\begin{lemma}
1) $r_C$ is a Lipschitz continuous, bounded  plurisubharmonic function on $\tM$. It is also an exhaustion function on $\tM$ if $d_C$ is complete.\\
2)  $r_C^2$ and $\log r_C^2$ are plurisubharmonic . \\
\label{rcproperty}
\end{lemma}

\begin{proof}
1) The boundedness of $r_C$ follows by its definition. 
To see that that $r_C$ is exhaustion, note that the Carath\'eodory distance function $d_{C,M}$ on $\widetilde{M}$ is complete by assumption.
Here $\ell_C$ actually approaches $+\infty$ and thus $r_C$ approaches $1$. It is well-known that $\ell_C$ is continuous (c.f. \cite[Proposition 3.1.13]{Koba1998}). Its Lipschitz continuity of $r_C$ follows from that of $\ell_C$, which may be found for example from \cite[\S 1.2]{Yeun2022}. 

Write $\ell=\ell_C, r=r_C$. To see that $r$ is plurisubharmonic, note that $d_C$ is obtained by taking supremum
among a set of plurisubharmonic functions, $\ell$ is a plurisubharmonic function on $\tM$. In fact, we have 

\[
\sqrt{-1}\ddbar \ell 
= \frac{1+\tanh^2\ell}{\tanh\ell} \sqrt{-1}\partial\ell\wedge \barD\ell
= \frac{1+r^2}{r}\sqrt{-1}\partial\ell\wedge \barD\ell 
\geq \sqrt{-1}\partial\ell\wedge \barD\ell,
\]
c.f. \cite{Yeun2022} or \cite{WY2022}). In above, we have used the fact that $\partial\ell\wedge \barD\ell$ exists as a current, cf. the first paragraph of \cite[\bf 1.3]{Yeun2022}.
Then by direct computation,
\begin{align*}
\sqrt{-1}\ddbar r_C 
= \sqrt{-1}\partial (\barD \tanh\ell)
&=\sqrt{-1}\partial(\sech^2\ell\barD\ell)  \\
&=-2\tanh\ell \sech^2\ell \sqrt{-1}\partial\ell\wedge \barD\ell+\sech^2\ell \sqrt{-1}\ddbar\ell \\
&= \bigg(-2r(1-r^2) +(1-r^2)\frac{1+r^2}{r}\bigg) \sqrt{-1}\partial\ell\wedge \barD\ell\\
&= \frac{1-r^4}{r} \sqrt{-1}\partial\ell\wedge \barD\ell \\
&\geqslant 0.
\end{align*}

2) Observe that at points where the following expressions are twice differentiable,
\begin{align*}
\barD r^2=\barD \tanh^2 \ell  
&= 2\tanh\ell\sech^2\ell  \barD\ell = 2r(1-r^2)\barD\ell \\ 
\sqrt{-1}\ddbar r^2
&=2(1-r^2)\left[(1-3r^2)\sqrt{-1}\partial\ell\wedge \barD\ell+r\sqrt{-1}\ddbar\ell \right] \\
&=2(1-r^2)\left[(1-3r^2)+(1+r^2) \right] \sqrt{-1}\partial\ell\wedge \barD\ell 
\tag{$\because \sqrt{-1}\ddbar\ell = \frac{1+r^2}{r}\sqrt{-1}\partial\ell\wedge \barD \ell$}\\
&=4(1-r^2)^2 \sqrt{-1}\partial\ell\wedge \barD\ell \\
&\geqslant 0. 
\end{align*}
So
\begin{align*}
\sqrt{-1}\ddbar \log r^2
&= \sqrt{-1}\partial (\frac{\barD r^2}{r^2})
= \frac{-1}{r^4}\sqrt{-1}\partial r^2\wedge \barD r^2 + \frac{1}{r^2}\sqrt{-1}\ddbar r^2 \\
&= \frac{-1}{r^4} \cdot 4r^2(1-r^2)^2 \sqrt{-1}\partial \ell \wedge \barD \ell 
+  \frac{1}{r^2} 4(1-r^2)^2 \sqrt{-1}\partial\ell\wedge \barD\ell \\
&=0.
\end{align*}

In general, the above arguments work for the Poincar\'e disk $\Delta$.   In  expressions such as $\pd\op \ell^2$, $\ell$ is taken as supremum of $f^*\ell_P$ for the corresponding Poincar\'e length function
$\ell_P$ on $\Delta$ for $f:M\rightarrow \Delta$ as defined earlier, and hence is plurisubharmonic.  The expression $\pd\op \ell^2$ is considered as a current.

\end{proof}

\begin{remark}
In \cite{Yeun2022}, it is also shown that $-\log(r^2-1)$ is plurisubharmonic.
\end{remark}

Let $a=a(R)= \tanh (cR)$.
Then 
\begin{equation}
B_{d_C}(x;cR)=\{y\in \tM \mid \ell_C(y)<cR\}=\{y\in \tM \mid r_C(y)< a\}.
\label{Bdc}
\end{equation}
Recall the notations in Proposition \ref{lvol}. It is reduced to the following lemma.

\begin{lemma}  There exists a constant $\alpha_0>0$ such that
$\Vol_{{KE}}(V\cap B_{_{KE}}(x;R))\geqslant \alpha_0\cdot \frac 1{1-a}=\alpha_0\cdot \frac{1}{1-\tanh (cR)}$, for $R>0$ sufficiently large.
\label{lemvolest}
\end{lemma}
The proof will be similar to  \cite[(2.3.13)]{HT2002} once we have the appropriate setting. We will use the following facts which may be found for example from \cite[Proposition 2.2.1]{HT2002}:
\begin{proposition}
Let $V\subset \tM$ be a $k$-dimensional complex analytic subvariety and $[V]$ be the closed $(n-k,n-k)$-current corresponding to $V$.
Denote by $\nu_x(\eta)$ the Lelong number at $x$ of a function $\eta$ which is plurisubharmonic on some relative compact open subset $W\subset \tM$ containing $x$; and by $\mult_x(\eta)$  the multiplicity of $\eta$ at $x$. Then \\ 
(i) \[
\int_{W} [V]\wedge (\sqrt{-1}\ddbar \eta)^k \geq \mult_{x}(V)\cdot \nu_x(\eta)^k.
\]\\
(ii) Let $\rho$ be function smooth on $W$ such that $\rho\equiv \eta$ on $W-W'$, where $W'\subset W$ is a relative compact open subset. Then
\[
\int_{W} [V]\wedge (\sqrt{-1}\ddbar \rho)^k = \int_{W} [V]\wedge (\sqrt{-1}\ddbar \eta)^k.
\]
\label{lelong}
\end{proposition}
\subsubsection*{Proof of Lemma \ref{lemvolest}  }
\begin{proof}
Let $\epsilon=1-a$.  Consider a smooth cut-off function $\chi:[0,1]\rightarrow[0,1]$ having the following properties: 
(i) $\chi$ is supported on $[0,1-\epsilon)$;
 (ii) $\chi(t)=1$ for $t\leqslant 1-2\epsilon$;
 (iii) $\chi$ is decreasing;
 (iv) $|\chi'|<\frac2{\epsilon}$; and
 (v) $|\chi^{\prime\prime}|\leqslant \frac2{\epsilon^2}$.
For $t\in (0,1]$,  consider the function 
\[
l(t):= \chi(t)\log t.  
\]
It follows that $l$ is increasing on $[0,1]$, i.e.,  $l'(t)\geq 0$ on $[0,1]$.
For $\frac{1}{2}>\epsilon>0$ sufficiently small and $t\in (1-2\epsilon,1-\epsilon)$,
 \begin{eqnarray*}
 |l'|
  & \leqslant& |\frac{\chi(t)}{t}| + |\chi'(t)||\log t|\\
  &\leqslant& \frac1{1-2\epsilon}+\frac 2\epsilon \cdot|\log(1-2\epsilon)|\leqslant   2+2\frac{|\log(1-2\epsilon)|}{\epsilon} \leqslant 6 \\
  |l^{\prime\prime}|
  &\leqslant & |\frac{\chi(t)}{t^2}| +2|\frac{\chi(t)}{t}|+|\chi^{\prime\prime}(t)||\log t|  \nonumber \\
  &\leqslant& \frac{1}{(1-2\epsilon)^2} + 2\cdot \frac{2}{\epsilon}\cdot \frac{1}{1-2\epsilon} + \frac{2}{\epsilon^2}|\log(1-2\epsilon)|  \nonumber \\
  &\leqslant & 4+  \frac{2}{\epsilon}\bigg(4+\frac{|\log(1-2\epsilon)|}{\epsilon} \bigg)
  \leqslant \frac{14}{\epsilon}.
  \end{eqnarray*}
Consider now 
\[
\psi:=l(r_C^2).
\]
For $r_C^2< 1-2\epsilon$, $\sqrt{-1}\ddbar \psi = \sqrt{-1}\ddbar \log r^2 \geqslant 0$ by Lemma \ref{rcproperty}. 

For $r_C^2\in [1-2\epsilon,1-\epsilon)$,
\begin{eqnarray}
\sqrt{-1}\pd\op\psi&=& \sqrt{-1}\pd\op l(r_C^2) \nonumber \\
&=& \sqrt{-1}\pd( l'(r_C^2) \cdot 2r_C\op r_C)\nonumber \\
&=&  l^{\prime\prime}(r_C^2) \cdot  4r_C^2\cdot \sqrt{-1}\pd r_C\wedge\op r_C
   + l'(r_C^2)\cdot 2\cdot \sqrt{-1}\pd r_C\wedge\op r_C
    +l'(r_C^2)\cdot 2 r_C \cdot \sqrt{-1}\pd\op r_C\nonumber \\
&\geqslant & 4r_C^2 \cdot l^{\prime\prime}(r_C^2) \cdot \sqrt{-1}\pd r_C\wedge\op r_C
\end{eqnarray}
where we have used the fact that $l$ is increasing and that $i\pd\op r_C$ is a positive current (Lemma \ref{rcproperty}). Then
  \[
  l^{\prime\prime}(r_C^2)\cdot 4r_C^2 
  \geq -\frac{14}{\epsilon}\cdot 4(1-\epsilon) 
  \geqslant -\frac{56}{\epsilon}.
  \]
Similar to the K\"ahler form of the Poincar\'e metric on $\Delta$, define 
\[
\tomega:=\frac{\sqrt{-1}\pd r_C\wedge \op r_C}{(1-r_C^2)^2}.
\]  
Note that for $r_C^2\in (1-2\epsilon,1-\epsilon)$,
\[
\frac{1}{(1-r_C^2)^2}=\frac{1}{(1+r_C)^2(1-r_C)^2}>\frac{1}{4}\cdot \frac{1}{4\epsilon^2}= \frac{1}{16\epsilon^2}.
\]
For $r_C\in (1-2\epsilon,1-\epsilon)$, let $\alpha>56\cdot 16$ be a constant, we get 
  \begin{align}
 \sqrt{-1}\pd\op(\frac{1}{\alpha \epsilon} \psi)+\tomega 
\geqslant & \left[ \frac{1}{\alpha \epsilon} l^{\prime\prime}(r_C^2)\cdot 4r_C^2 +\frac{1}{(1-r_C^2)^2} \right] \sqrt{-1}\pd r_C\wedge\op r_C \nonumber \\
>& \left[ \frac{1}{\alpha\epsilon} (-\frac{56}{\epsilon}) +\frac{1}{16\epsilon^2} \right] \sqrt{-1}\pd r_C\wedge\op r_C>0.
\label{psiest} 
  \end{align}

By definition, for $x\in \tM$ and $v\in T_x\tM$,
\[
  g_C(x;v)\geqslant\sup_{f\in \mathcal{F}} \{ f^*g_{\Delta,P}(x;v)\},
\] 
where we have denoted by $\mathcal{F}$ the family of holomorphic functions $f:\tM\rightarrow \Delta$ with $f(x)=0$; and by $g_{\Delta,P}$ the Poincar\'e metric on $\Delta$. 

Write $\omega_{_{KE}}$ as the K\"ahler form associated to $g_{_{KE}}$.
Together with \eqref{CKEest}, we get
\begin{equation}
\omega_{KE}
\geqslant c\cdot \sup_{f\in \mathcal{F}}f^*\omega_{\Delta,P}
=c\cdot\sup_{f\in \mathcal{F}}f^*\left(\frac{|dz|^2}{(1-|z|^2)^2}\right)
=c\cdot\frac{|d r_C|^2}{(1-r_C^2)^2}
=c\cdot\tomega.
\label{omegaKEest}
\end{equation}
Define on $\tM$ the function
\[
\phi:=\frac{1}{c\alpha (1-a)}\psi,
\]
which is compactly supported in $\{r_C<1-\epsilon=a\}\subset \tM$.
It follows from \eqref{psiest} and \eqref{omegaKEest} that
\[
i\pd\op\phi+\tomega \geqslant 0\quad \text{on $\tM$}.
\]
On $\tM$, which is simply connected, we may write 
\begin{equation}
\omega_{KE}=\sqrt{-1}\pd\op\Phi 
\label{KEpotential}
\end{equation}
in terms of a potential function $\Phi$.
  Then we conclude that 
\[
  \phi+\Phi
\]
is a plurisubharmonic function on $\tM$.

Now note that 
\begin{equation}
\nu_x(\phi)=\frac1{c \cdot\alpha (1-a)}.
\label{lelongnu}
\end{equation}
Then
\begin{align*}
&\Vol_{_{KE}}(B_{_{KE}}(x;R)\cap V) \\
=& \int_{B_{_{KE}}(x;R)} [V]\wedge \omega_{_{KE}}	\tag{by definition}\\
\geqslant&  \int_{B_{d_C}(x;cR)} [V]\wedge \omega_{_{KE}}  \tag{by \eqref{Bsubset}} \\
=&  \int_{\{r_C<a\}} [V]\wedge \omega_{_{KE}}  \tag{by \eqref{Bdc}} \\
=& \int_{\{r_C<a\}} [V]\wedge (\sqrt{-1}\ddbar \Phi)  \tag{by \eqref{KEpotential}} \\
=& \int_{\{r_C<a\}} [V]\wedge (\sqrt{-1}\ddbar (\phi+\Phi)) \tag{by Proposition \ref{lelong} (ii), since $\phi$ has compact support} \\
\geqslant & \mult_x(V) \cdot \nu_x(\phi+\Phi) \tag{by Proposition \ref{lelong} (i) and the plurisubharmonicity of $\phi+\Phi$} \\
=& \mult_x(V) \cdot \bigg( \nu_x(\phi) +\nu_x(\Phi)\bigg) \tag{by the definition of Lelong number} \\
=& \mult_x(V) \cdot \frac{1}{c\alpha(1-a)} \tag{$\nu_x(\Phi)=0$ since $\Phi$ is smooth at $x$},
\end{align*}
where $\alpha>56\cdot 16$.
\end{proof}

\begin{remark}
By using \cite[Proposition 3.1.2]{HT2000}, it is possible to obtain a sharper volume estimate in Lemma \ref{lemvolest}.  The above gives a more direct
construction sufficient for our purpose.  In the Appendix, we will give another possible construction for $\phi$ given by $\phi=\psi_\varepsilon(r_C)$, where $\psi_\varepsilon$ is as in  \cite[Proposition 3.1.2]{HT2000}.
\end{remark}

\begin{proof}[Proof of Proposition \ref{lvol}]
It is an immediate consequence of Lemma \ref{lemvolest}.
\end{proof}

\begin{proof}[Proof of Theorem \ref{MT3}]
This now follows from Proposition \ref{lvol}.
\end{proof}

\subsection{Upper bound in volume}
To obtain an upper bound for the volume of curves, 
an important step as in \cite[Proposition 3.1]{HT2006} is to apply Royden's Schwarz Lemma \cite{Royden1980} to get a comparison between the pullback canonical metric and the hyperbolic metric on a Riemann surface. In our situation, we will need to replace Royden's Schwarz Lemma by a more general Schwarz Lemma applicable to nonsmooth complex Finsler metrics:
\begin{lemma}
Let $X$ be a complex manifold equipped with a nondegenerate infinitesimal Carath\'eodory metric $g_C$. Let $R=\Delta/\Gamma$ be a compact Riemann surface of genus $\geqslant 2$, equipped with the Hermitian metric $g_R$ obtained by descending the Poincar\'e metric  on $\Delta$. 
Suppose $\phi: (R,g_R) \rightarrow (X,g_C)$ is a nonconstant holomorphic map. Assume that the following conditions are satisfied: \\
i) the Gauss curvature of $g_R$ is bounded from below by $-k_R$ for some $k_R>0$; \\
ii) the holomorphic sectional curvature of $g_C$ is bounded from above by $-k_C$ for some $k_C>0$.\\
Then $\phi^*g_C\leqslant \frac{k_R}{k_C} g_R$.
\label{schwarz}
\end{lemma}

\begin{proof}
Let $w$ be the local coordinate on $R$. Write $\phi^*g_C:=d\sigma^2= 2 \lambda dw\otimes d\overline{w}$. For the K\"ahler metric $g_R$, we write $g_R= 2\mu dw\otimes d\overline{w}$. Both $\lambda$ and $\mu$ are nonnegative. Let $u = \frac{\phi^*g_C}{g_R}= \frac{\lambda}{\mu}$. It suffices to show that $u\leqslant \frac{k_R}{k_C}$ on $R$.

The infinitesimal Carath\'eodory metric $g_C$ is upper-semicontinuous, so the pullback $\phi^*g_C:=d\sigma^2$ is an upper-semicontinuous Hermitian pseudo-metric on the compact Riemann surface $R=\Delta/\Gamma$. 
Thus $\lambda$ is upper-semicontinuous on $R$ and so is $u = \frac{\lambda}{\mu}$. 
Since $R$ is compact, there exists $w_0\in R$ such that $u$ attains its maximum at $w_0$.
Since $g_C$ is plurisubharmonic (thus so is $\lambda$ and $u$), 
we can take the Laplacian of $u$ in the sense of a current, or a distribution.
Hence as $M$ has complex dimension $1$, we understand in the following that for functions $f$ relevant to our discussions,
$$\frac{\partial^2 f}{\partial w \partial\bar w}(w_0):=\frac14\liminf_{r\rightarrow 0}\frac1{r^2}\int_{0}^{2\pi}(f(w_0+re^{i\theta})-f(w_0))\frac{d\theta}{2\pi}.$$
With this interpretation, as $w_0$ is also a maximum point of $u$,   the Maximum Principle implies that
\[
0\geqslant \frac{\partial^2 \log u}{\partial w \partial\bar w}(w_0) 
=\frac{\partial^2 \log \lambda}{\partial w \partial\bar w}(w_0) -\frac{\partial^2 \log \mu}{\partial w \partial\bar w}(w_0)
=-\lambda(w_0) K_{d\sigma^2}(w_0)+\mu (w_0) K_{g_{R}}(w_0).
\]
Note that the holomorphic sectional curvatures
\[
K_{d\sigma^2}(v)=K_{\phi^*g_C}(v)\leqslant K_{g_C}(d\phi(v))\leqslant -k_C<0, \quad v\in TR,
\]
see \cite[p. 31-32]{Koba1998} or \cite{WW2011}). It follows that
\[
u(w_0)=\frac{\lambda(w_0)}{\mu (w_0)}
\leqslant \frac{K_{g_R}(w_0)}{K_{d\sigma^2}(w_0)}
\leqslant \frac{k_R}{k_C}.
\]
Since $w_0$ is a maximum point of $u$ on $S$, it follows that $u\leqslant \frac{k_R}{k_C}$ on $R$.
\end{proof}

We give an alternate argument using the technique of Ahlfors on the proof of Schwarz Lemma \cite{Ahl1938}.
\begin{proof}[Alternative proof of Lemma \ref{schwarz}]
Let $\pi: \widetilde{R}\rightarrow R$ be the universal covering map.
Denote by $\Phi: \widetilde{R}\cong \Delta \rightarrow X$ the lifting of $\phi:R\rightarrow X$.
We also let $g_{\tR}$ to denote the Poincar\'e metric on $\widetilde{R}$. 
\[
\begin{tikzcd}
(\widetilde{R},g_{\tR})  \arrow[d,"\pi"] \arrow[dr,"\Phi"] & \\
(R,g_R) \arrow[r,"\phi"] & (X,g_C)
\end{tikzcd}
\]
Suppose $w_0$ is a point on $R$ where $u=\frac{\phi^*g_C}{g_R}$ has maximal value.
Let $z_0=\phi(w_0)$.  From definition of Carath\'eodory metric and a normal family argument, there exists a mapping $h_{z_0}: X\rightarrow\Delta$ such
that $h_{z_0}(z_0)=0$ and $g_C(z_0)=h_{z_0}^*g_{P,\Delta}(z_0)$.
Here $g_{P,\Delta}$ denotes the Poincar\'e metric on $\Delta$.

Let $\tw_0\in \tR$ be a point so that $\pi(\tw_0)=w_0$.
Note that $g_{\tR}=\pi^*g_R$, so 
\[
K_{g_{\tR}}(\tw)=K_{g_R}(\pi(\tw))\geqslant -k_{R}, \qquad  \forall \tw\in \tR.
\]
Consider now the function 
\[
\tu_{z_0}(\tw):=\frac{\pi^*\phi^*h_{z_0}^*g_{P,\Delta}}{\pi^*g_R}(\tw)
\]
on $\tR$.  
If $\tu_{z_0}$ achieves a maximum at a point $\tw_1\in \tR$, 
applying $\ddbar$ and arguing using Maximum Principle as in previous proof shows that 
$\tu_{z_0}(\tw)\leqslant \tu_{z_0}(\tw_1)\leqslant \frac{k_R}{k_C}$ all $\tw\in \tR$.
This implies in particular that for all $w\in R$,
\begin{align*}
u(w)\leqslant u(w_0)
=\frac{\phi^*g_C(w_0)}{g_R(w_0)} 
= \frac{\phi^*g_C(\pi(\tw_0)}{g_R(\pi(\tw_0))}
=\frac{\pi^*\phi^*h_{z_0}^*g_{P,\Delta}}{\pi^*g_R}(\tw_0)
=\tu_{z_0}(\tw_0)
\leqslant \frac{k_R}{k_C}.
\end{align*}

In general, to find a maximum point, we apply the trick of using barrier as Ahlfors \cite{Ahl1938}.  For $1>a>0$, the Poincar\'e metric on $\Delta_a:=\{z\in \bC:|z|<a\}$ is given
by 
\[
g_{\Delta_a}(w)=\frac{|d(\frac wa)|^2}{(1-|\frac wa|^2)^2}=\frac{a^2|d w|^2}{(a^2-|w|^2)^2}.
\] 
Note that in the above discussion, $\pi^*g_R=g_{\Delta}$.  
Instead of 
$\tu_{z_0}(\tw)=\frac{\pi^*\phi^*h_{z_0}^*g_{P,\Delta}}{\pi^*g_R}(\tw)$, we consider 
\[
\tu_{a,z_0}(\tw):=\frac{\pi^*\phi^*h_{z_0}^*g_{P,\Delta}}{g_{\Delta_a}}(\tw).
\]  
As $|\tw|\rightarrow a$, $\pi^*\phi^*h_{z_0}^*g_{P,\Delta}(\tw)$ is  bounded while 
$g_{\Delta_a}(\tw) \rightarrow \infty$.
We see that the supremum of $\tu_{a,z_0}$ has to be achieved at a point $\tw_{a}$ lying in the interior of $\Delta_a$.    Hence the above argument implies that for any $\tw\in \Delta_a$,
\[
\tu_{a,z_0}(\tw)\leqslant \tu_{a,z_0}(\tw_a)\leqslant \frac{k_R}{k_C},
\]
where the right hand side is independent of $a$.  Letting $a\rightarrow 1$, we conclude that $\tu(\tw)\leqslant k$ for all 
$\tw\in \Delta\cong \tR$.  The rest of the  argument is the same as before.

\end{proof}

Recall the notations of Theorem \ref{MT1}. Write $M=M_i$ for simplicity.
We may now state the following slight modification of \cite[Proposition 3.1]{HT2006}:
\begin{proposition}
Let $S=\Delta/\Gamma$ be a compact Riemann surface of genus $g(S)\geqslant 2$.
Suppose $f:S\rightarrow \overline{M}$ is a nonconstant holomorphic map such that $f(S)\cap M\neq\emptyset$. Let $w\in S$ such that $x=f(w)\in f(S)\cap M$. Then 
there exists a constant $k>0$ such that
for any $R>0$,
\begin{eqnarray*}
 \Vol_{{_{KE}}}\bigg(f(S)\cap B_{_{KE}}(x;R)\bigg)
\leqslant k (2g(S)-2).
\end{eqnarray*}
\label{uvol}
\end{proposition}

\begin{proof}
Equip $S=\Delta/\Gamma$ with a K\"ahler metric $h$, which is obtained by descending the 
Poincar\'e metric on $\Delta$.
We may suppose $h$ is of constant Gaussian curvature $-k_h$ for some $k_h>0$.
Since $g_C\leqslant -k_C$ for some constant $k_C>0$ \cite{Bur1978}, we may apply Lemma \ref{schwarz} to  $f_1:=f|_{f^{-1}(M)}: (f^{-1}(M),h)\rightarrow (M,g_C)$ to see that $f_1^*g_C\leqslant \frac{k_h}{k_C}h$. By the assumption \eqref{CKEest} that
$g_{_{KE}}\leqslant \frac{1}{c}\cdot g_C$ for some $c>0$ on the universal covering $\tM$ of $M$.  Since both $g_{_{KE}}$ and $g_C$ are invariant under $\Aut(\tM)$, the last inequality descends to hold on $M$.
But then $f_1(p)=f(p)$ for any $p\in f^{-1}(M)$. Thus 
\[
f^*g_{_{KE}}\leqslant \frac{1}{c} f^*g_C \leqslant \frac{k_h}{k_C}h, \quad \text{on $f^{-1}(M)$}. 
\]
Letting $k=c\cdot \frac{k_h}{k_C}$, we have
\begin{eqnarray*}
 \Vol_{{_{KE}}}\bigg(f(S)\cap B_{_{KE}}(x;R)\bigg)
&=& \int_{f(S)\cap B_{_{KE}}(x;R)} \omega_{_{KE}} \\
&\leqslant&  \int_{f^{-1}(M)} f^*\omega_{_{KE}} 
\leqslant  \int_{f^{-1}(M)} k \omega_h 
\leqslant k \int_S  \omega_h 
=k (2g(S)-2).
\end{eqnarray*}
\end{proof}

\begin{remark}
In the case that $M$ is a Hermitian locally symmetric space, $g_{_{KE}}$ is the canonical metric induced by the Bergman metric on $\tM$, whose holomorphic sectional curvature is bounded from above by some constant $-k_M$ for some $k_M>0$. One may apply 
Royden's Schwarz Lemma \cite{Royden1980} to $f_1:=f|_{f^{-1}(M)}: (f^{-1}(M),h)\rightarrow (M,g_{_{KE}})$ to conclude that $f_1^*g_{_{KE}} \leqslant \frac{k_M}{k_h} h$. In our case, the required negative holomorphic sectional curvature upper bound on $(M,g_{_{KE}})$ is not clear.
\end{remark}

Write $V=f(S)$.
As is discussed at the beginning of this section, for finding volume bounds, it is equivalent to consider 
 local liftings $\widetilde{V}$ of $V\cap M$ 
to the universal covering $\tM$ of $M$. 
We would use the notation $V=\widetilde{V}$ on $\tM$ for the such local liftings for simplicity.
Combining Proposition \ref{lvol} (or Lemma \ref{lemvolest}) and Proposition \ref{uvol}, we have on $\tM$ and for $R>0$ sufficiently large,
\begin{equation*}
(*):\quad 
\frac{\alpha_0}{1-\tanh(c\cdot R)}
\leqslant \Vol_{{_{KE}}}\bigg(V\cap B_{_{KE}}(x;R)\bigg)
\leqslant \frac{B}{A}(2g(S)-2).
\label{*}
\end{equation*}

\section{ Injectivity radius lower bound}
\label{inj}
Suppose $f(S)\cap M_i\neq \emptyset$.
Denote by $\rho_i(x)$ the injectivity radius at $x\in f(S)\cap M_i$ with respect to $g_{_{KE}}$. If \hyperref[*]{$(* )$} hold for $R=\rho_i(x)$, then there exist fixed real numbers $\alpha,\beta,\gamma$ such that 
\begin{equation}
\rho_i(x) \leq \alpha\log(\beta g(S)+\gamma)=:\tau_0.
\label{injbdd}
\end{equation}
This implies that if $w\in S$ with $f(w)=x$ is such that $\rho_i(x)>\tau_0$, then $x\not\in f(S)\cap M_i$. For proving Theorem \ref{MT1}, it suffices to use this observation on sufficiently large coverings $M_i\rightarrow M$. We first need the following:
\begin{proposition}
Fix $g\geqslant 2$. There exists a compact subset $Y\subset M_i=\tM/\Gamma_i$ having the following property:  for any subgroup $\Gamma_i\triangleleft \Gamma$ of finite index, if there is a compact Riemann surface $S$ of genus $g$ and a nonconstant holomorphic map $f: S\rightarrow M_i$, then there is a compact Riemann surface $S'$ of genus $g'\leqslant g$ and a nonconstant holomorphic map $f':R'\rightarrow M_i$ such that $\pi_i\circ f'(S')\cap Y\neq \emptyset$. Here $\pi_i:M_i\rightarrow M$ is the covering map.
\label{cptset}
\end{proposition}
\begin{proof}
By assumption, $M=\tM/\Gamma$ is quasi-projective variety with ample canonical line bundle. By Proposition \ref{uvol}, $V=f(R)$ is a projective curve in $M$ of degree bounded by $C_1\cdot g(R)$ for some $C_1>0$. Then the deformation theoretic argument in the proof of \cite[Proposition 3.4]{HT2006} may be applied.
\end{proof}

Following \cite[\S2]{HT2006},
a set of subgroups $\{\Gamma_i\leq \Gamma: i\in \Lambda\}$ in $\Gamma$ is said to be separating if for each infinite subset $J\subset \Lambda$, $\bigcap_{j\in J} \Gamma_i=\{1\}$. Since $M$ supports a tower of covering $\{M_i\}_{i\in I}$, it follows that $\Gamma$ has a separating set of subgroups indexed by $I$.

\begin{proposition}
Let $S$ be a compact Riemann surface of genus $g$.
Suppose $f: S\rightarrow \overline{M_i}$ is a nonconstant holomorphic map . Then there exists a compact subset $Z\subset M_i$ and $i_0\geqslant 0$ such that for all $i\geqslant i_0$, the injectivity radius
$\rho_{i}(x) >\tau_0$ for any $x\in Z$, and there exists a curve $S'$ of genus $g'\leqslant g$ and $f:S'\rightarrow M_i$ such that $\pi_i\circ f'(S')\cap Z\neq\emptyset$.
\label{injlowbdd}
\end{proposition}

\begin{proof}
We have the finite unramified covering $\pi_i:M_i\rightarrow M$. By Proposition \ref{cptset}, there is a compact Riemann surface $S'$ with genus $g(S')\leqslant g(S)$ and a nonconstant holomorphic map $f': S'\rightarrow M_i$ whose image $\pi_i\circ f'(S')$ always cut a compact subset $Y\subset M$. Let $Z\subset M_i$ be a compact set so that $Y\subset \pi_i(Z)$. Since $\{\Gamma_i\}_{i\in I}$ forms a set separating subgroups of $\Gamma$, by \cite[Proposition 2.3]{HT2006}, there exists a finite subset $I_{Z}\subset I$ such that for any $x\in Z$ and any $j\in I\backslash I_{Z}$, we have $\rho_j(x)>\tau_0$. For $i_0\geqslant 0$ sufficiently large, we can make sure $\{\Gamma_i: i\geqslant i_0\}\cap I_Z =\emptyset$ for all $i\geqslant i_0$ since $\bigcap_{i\in I}\Gamma_i=\{1\}$. Hence the result follows.
\end{proof}

\section{Ramification index lower bound}
\label{ramind}
In this section, we are going to generalize \cite[Proposition 4.4]{HT2006}. For $\bar\pi: \overline{M_i}\rightarrow \overline{M}$ extending $\pi: M_i\rightarrow M$, write $D_i=\overline{M_i}-M_i$ and $D=\overline{M}-M$. We will follow the strategy of \cite[\S 4]{HT2006}. In the following, we suppose $i\geqslant 0$ is sufficiently large.

\begin{proposition}
For any $x\in M_i$, there exists $q_0\geqslant 0$ such that whenever $q\geqslant q_0$,
we have a section $s\in H^0(\overline{M_i},q(K_{\overline{M_i}}+D_i))$ so that $s(x)\neq 0$ and $s|_{D_i}\equiv 0$.
\label{sec}
\end{proposition}
\begin{proof}
By \cite{WY2022}, $\overline{M}$ is of log-general type with respect to $D$, i.e. $K_{\overline{M}}+D$ is big. In fact, there exists $q_0\geqslant 0$ such that for $q\geqslant q_0$, we have a nontrivial section $\sigma\in H^0(\overline{M},q(K_{\overline{M}}+D))$ where the order of jets of $\sigma$ at $\bar\pi_i(x)$ can be prescribed up to order $cq^n$ for some constant $c>0$ (c.f. proof of \cite[Theorem 0.6]{WY2022}).
Write $D_i=\cup_{j} D_i^j$ as the union of irreducible components $D_i^j$'s. 
Let $r^j$ be the ramification index of $\bar\pi_i$ at $D_i^j$'s and $m:=\min_j r^j$.
Note that
\begin{align*}
\bar\pi_i^* \bigg(q(K_{\overline{M}}+D) \bigg)
&= q (K_{\overline{M}}+D_i)- \sum_j(r^j-1)D_i^j \\
&\leq q (K_{\overline{M}}+D_i)- (m-1)D_i
\end{align*}
Therefore $s:=\bar\pi_i^* \sigma \in H^0(\overline{M_i}, q (K_{\overline{M_i}}+D_i)- (m-1)D_i)$. Choose an $i\geqslant 0$ sufficiently large so that $m>q+1$. Then the above implies that $s$ vanishes on $D_i$. 
\end{proof}

For the quasi-projective manifold $M=\oM - D$, there exists a complete Poincar\'e-type K\"ahler metric $g_P$ of bounded geometry on $M$ whose construction is recalled as follows. 
In a neighbourhood $U\cong \Delta^n$ of a point on $D$ in $\oM$, there is the Poincar\'e metric $\rho$ on $\Delta^n$ whose associated K\"ahler form is
\begin{equation}
\eta:=\frac{\sqrt{-1}dz_1\wedge d\overline{z_1}}{|z_1|^2(\log|z_1|^2)^2}+\dots+\frac{\sqrt{-1}dz_k\wedge d\overline{z_k}}{|z_k|^2(\log|z_k|^2)^2}+ \sqrt{-1}dz_{k+1}\wedge d\overline{z}_{k+1}+\dots +\sqrt{-1}dz_{n}\wedge d\overline{z}_n. 
\label{poin}
\end{equation}
for $|z_i|<\frac12$.
The Poincar\'e type metric $g_P$ on $\overline{M}$ is obtained by patching up the above metrics $\rho$'s on a finite number of such neighborhoods of $D$ together with
a smooth metric on the complement of the union of the neighborhoods on $M$, by partition of unity. In the following, we also denote  by $\omega_P$ the corresponding K\"ahler form of $g_P$ on $\overline{M}$.

\begin{proposition}
Let $f:S\rightarrow \overline{M_i}$ be a nonconstant holomorphic map such that $f(S)\cap M_i\neq \emptyset$. Then there is a constant $A_0>0$ depending on the holomorphic sectional curvature of $g_{_{KE}}$, such that
\[
\deg(f^* c_1(K_{\overline{M_i}})) + \deg (f^*c_1(D_i)) \leq A_0\cdot \bigg[\deg(c_1(K_S))+\deg (f^*c_1(D_i))\bigg]
\]
\label{deg}
\end{proposition}
\begin{proof}
The idea of proof is standard by now and can be found in \cite{Nade1989}  and \cite{Nogu1991}.  We give outline here in our setting, for completeness of presentation.

By the assumption in Theorem \ref{MT1}, $M_i$ is equipped with a K\"ahler-Einstein metric $g_{_{KE}}$. Denote by $\omega_{_{KE}}$ the associated K\"ahler form of $g_{_{KE}}$, normalized so that
\[
\frac{\sqrt{-1}}{2\pi}\ddbar \log \det g_{_{KE}} = \omega_{_{KE}}.
\]
Let $g_P$ be the Poincar\'e-type metric on $\overline{M_i}$ with K\"ahler form $\omega_P$. We recall some facts about $g_P$ which can be found for example from \cite{WY2022}.
By Royden's Schwarz Lemma \cite{Royden1980}, $g_{_{KE}}\leqslant c'g_P$ for some constant $c'>0$.

 In view of \eqref{poin}, $g_P$ has pole along $D$ of order $\leqslant 2$ and so is $g_{_{KE}}$. 
Let $\sigma$ be a local holomorphic section vanishes along $D_i$.
Let $g_o$ be a smooth K\"ahler metric on $\overline{M}_i$. 
There exist a constant $c_0>0$ such that 
\[
\frac{\det g_o}{|\sigma|^2\det g_{_{KE}}}\leqslant c_0.
\]
So
\begin{align*}
{\frac{\sqrt{-1}}{2\pi}}\ddbar \log \bigg(\frac{\det g_o}{|\sigma|^2\det g_{_{KE}}} \bigg)
\geqslant 
f^*c_1(K_{\overline{M_i}}) +f^*c_1(D_i)-f^*\omega_{_{KE}}.
\end{align*}
as current on $S$. Thus
\begin{equation}
\deg (f^*\omega_{_{KE}})\geqslant \deg (f^*c_1(K_{\overline{M_i}}))+\deg (f^*c_1(D_i)).
\label{deg1}
\end{equation}

Let $\Phi$ be a volume  on $S$, which for instance could be taken as the volume  corresponding to the metric induced by the Poincar\'e metric on $\Delta$.
Then $c_1(K_S)=\frac{\sqrt{-1}}{2\pi}\ddbar \log \Phi$.
Note that there is $c_0'>0$ such that
\[
\frac{f^*(|\sigma|^2\omega_{_{KE}})}{\Phi} \leqslant c_0' \quad \text{on S,}
\]
again from the pole order estimate of $\omega_{KE}$ near the compactifying divisor.
Let
\[
\mathcal{R}:= f^{-1}(D_i)\cup \{z\in S: df_z=0\}.
\]
Since holomorphic sectional curvature is bounded from above by a negative constant $-\gamma$, we have 
\[
{\frac{\sqrt{-1}}{2\pi}}\ddbar \log \frac{f^* \omega_{_{KE}}}{dz\wedge d\overline{z}}
\geqslant A_1 f^*\omega_{_{KE}} \quad \text{on $S-\mathcal{R}$.}
\]
Here $A_1=c'\gamma$ for some constant $c'>0$.
Then
\[
\frac{\sqrt{-1}}{2\pi} \ddbar \log \frac{f^* (|\sigma|^2\omega_{_{KE}})}{\Phi}
\geqslant -f^*c_1(D_i)+A_1 f^*\omega_{_{KE}}-c_1(K_S),
\]
which implies that
\begin{equation}
\deg(f^*\omega_{_{KE}})\leqslant \frac{1}{A_1} \bigg(\deg(f^*c_1(D_i))+\deg(c_1(K_S)) \bigg)
\label{deg2}
\end{equation}
Combining \eqref{deg1} and \eqref{deg2} and take $A_0=\frac{1}{A_1}$, we are done.
\end{proof}
\bigskip

Let $m$ be the smallest number among the ramification indices along $D_i$'s. 
\begin{proposition}
Let $g\geqslant 0$ be fixed. Suppose $m>q_0A_0(2g-1)$, where $q_0$ and $A_0$ are the constants in Proposition \ref{sec} and \ref{deg}.
Suppose $S$ is a compact Riemann surface  genus $g(S)=g$ and $f:S\rightarrow \overline{M_i}$ is a nonconstant holomorphic map, 
such that $f(S)\cap  D_i\neq \emptyset$. Then  $f(S)\subset D_i$.
\label{ram}
\end{proposition}
\begin{proof}
Suppose $f(S)\cap M_i\neq \emptyset$, we can find a point $x\in f(S)\cap M_i$. Then $\pi(x)\in M$. By Proposition \ref{sec},  there exists a section $s\in H^0\bigg(\overline{M},q_0(K_{\overline{M}}+D)\bigg)$ such that $s(\pi(x))\neq 0$ and $s|_{D}\equiv 0$. So the pullback $\xi:=\bar\pi^*s\in H^0\bigg(\overline{M_i},q_0(K_{\overline{M_i}}+D_i)-mD_i\bigg)$. Since $s(\pi(x))\neq 0$, we have $\xi|_{f(S)}\not\equiv 0$ so that $\deg(f^*\xi )\geq 0$. It follows that
\begin{align*}
0
&\leqslant q_0\bigg(\deg(f^*c_1(K_{\overline{M_i}}) + \deg(f^*c_1(D_i)))\bigg) -  m\deg(f^*c_1(D_i)
 \\
&\leqslant q_0 A_0\cdot \bigg[\deg(c_1(K_S))+\deg (f^*c_1(D_i))\bigg] -  m\deg(f^*c_1(D_i))\tag{Proposition \ref{deg}}\\
&= q_0 A_0\cdot \bigg[2g-2+\deg (f^*c_1(D))\bigg] -  m\deg(f^*c_1(D_i)) \\
&= q_0 A_0 (2g-2) + (q_0A_0-m) \deg(f^*c_1(D_i)) 
\end{align*}
So the above implies
\begin{align*}
2g-2
&\geqslant \bigg(\frac{m}{q_0A_0}-1\bigg)\deg(f^*c_1(D_i)) \\
&> \bigg(2g-2 \bigg)\deg(f^*c_1(D_i)) 	\tag{since by assumption $m>q_0A_0(2g-1)$}
\end{align*}
Since $f(S)\cap D_i\neq \emptyset$, we have $\deg(f^*c_1(D_i))\geq 1$. Therefore the above inequality leads to a contradiction.
\end{proof}

\section{Proof of Main Theorems}
\label{pfMT}

\subsubsection*{Proof of Theorem \ref{MT1}} 
\begin{proof}
The proof of Theorem \ref{MT1} follows from the same strategy as \cite{HT2006}.
We consider in the following $g_0\geqslant 2$. 

By \cite[Lemma 4.3]{WY2023}, for $i$ sufficient large, we may assume that the ramification index $m$ is as large as we want. Then by Proposition \ref{ram}, it suffices to show that $f(S)\cap D_i\neq\emptyset$. Suppose this is not the case, i.e. we have $f(S)\subset M_i$. 
Because $f:S\rightarrow M_i$ is nonconstant,
by Proposition \ref{cptset}, we can find a compact Riemann surface $S'$ of 
 genus $g'\leqslant g_0$ and a nonconstant holomorphic map $f':S'\rightarrow M_i$ so that $f(S')\cap Z\subset M_i$ for some compact subset $Z\in M_i$. Now Proposition \ref{injlowbdd} implies that $\rho_i(x)>\tau_0$ for any $x\in Z$. Here $\tau_0$ is defined by the equality in \ref{injbdd}. But by Proposition \ref{lvol} and \ref{uvol} (or by \hyperref[*]{$(* )$}), $\rho_i(x)>\tau_0$ cannot be satisfied for any $x\in f(S')$. Hence we reached a contradiction.
\end{proof}

\subsubsection*{Proof of Theorem \ref{MT2}}
\begin{proof}
1).
Note that (a) is already given 
in \cite{HT2006}.
In our more general situation, 
we give a slightly different arguments for both (a) and (b) using Theorem \ref{MT1}.

Given a genus $g_0\geqslant 2$ curve $S$ and a holomorphic map $f:S\rightarrow \overline{M^k}$.
In general $\overline{M^k}$ is singular and Theorem \ref{MT1} is not directly applicable. Consider a smooth toroidal compactification $X^k\supset M^k$. It is well-known that $\overline{M^k}$ is a minimal compactification (c.f. \cite{Mok2012} for the nonarithmetic case). There exists a unique holomorphic map $\sigma_k: X^k\rightarrow \overline{M^k}$ restricting to the identity on $M^k$. The holomorphic map $f$ lifts to a holomorphic map $\bar f: S\rightarrow X^k$ such that $f=\sigma_k\circ \bar f$. Now since $X^k$ is smooth, we may apply Theorem \ref{MT1} to $\bar{f}$ and suppose $k$ is sufficiently large to conclude that $\bar f(S)\subset X^k - M^k$. Therefore $f(S)=\sigma_k\circ \bar f(S)\subset \overline{M^k}-M^k =:D_k$. Note that the boundary $D_k$ is stratified 
by disjoint unions of Hermitian locally symmetric spaces of strictly lower dimensions.
So $f(S)$ lies on exactly one such stratum. The same argument above may be applied repeatedly to conclude that $f$ must in fact be a constant map.

\medskip

2).
The argument for $M_g$ is similar but has some essential differences. 
 Statement (I) follows from Theorem \ref{MT1}. Focuses will therefore be put on (II).

The Deligne-Mumford compactification $\oM_g\supset M_g$ is obtained by
adding $\cup_{i=0}^{[\frac g2]}E_i$, where each $E_i$ is a divisor on $\oM_g$. 
Here we write $[\frac g2]$ as the integral part of $\frac{g}{2}$.
By \cite[p. 304]{Mum1983},
\begin{align}
E_0 &=\oM_{g-1,2}, \nonumber\\ 
E_i &\cong \oM_{i,1}\times \oM_{g-i,1}, \quad \qquad \forall \ 1\leqslant i<\frac g2, \nonumber\\
E_{\frac g2} &\cong(\oM_{\frac g2,1}\times \oM_{\frac g2,1})/\mathbb{Z}_2, \quad {\mbox{if}}  \ \ 2|g,
\label{DMbdy}
\end{align}
where $\mathbb{Z}_2=\mathbb{Z}/2\mathbb{Z}$ acts by permuting the two punctures.
 Here $\oM_{g,n}$ is the moduli space parametrizing curves with unordered marked points. For each level $m\in \bN$, we write similarly $D_m = \overline{M_{g,n}^m}-M_{g,n}^m = \pd M_{g,n}^m =\cup_{i=0}^{[\frac{g}{2}]}E_{i}^m$. 

We make the following observation: 
\begin{claim}
Let  $F:S\rightarrow \overline{M_{g,n}^m}$ be a nonconstant holomorphic map and $\pi_n:\overline{M_{g,n}^m}\rightarrow \overline{M_{g}^m}$ be the projection given by forgetting the punctures. Consider $F_n:=\pi_n\circ F:S\rightarrow \overline{M_g^m}$.  Then  if $F_n$ is nonconstant, $F(S)$ has to lie in  $\pd M_{g,n}^m$  for $m$ sufficiently large.
\end{claim}

\begin{cproof}
Consider the induced Torelli map: $j^m_g: M^m_g\rightarrow \mathcal{A}^m_g$, where $\mathcal{A}^m_g$ is the Siegel modular variety with canonical level $k$ structure.  The mapping $j^m_g$ extends to $j^m_g:\overline{M^m_g}\rightarrow \overline{\mathcal{A}^m_g}$, where $\overline{\mathcal{A}^m_g}\supset \mathcal{A}^m_g$ is a toroidal compactification. 
Now by Theorem \ref{MT1}, $j^m_g\circ F_n(S)\subset \overline{\mathcal{A}^m_g}-\mathcal{A}^m_g$. Apply case 1) of Theorem \ref{MT2}, i.e. the case of Hermitian locally symmetric spaces, we see that $j^m_g\circ F_n$ is a constant map for $m$ sufficiently large. Therefore $F_n(S)$ lies on a fibre $(j^m_g)^{-1}(b)\subset \overline{M^m_g}-M^m_g = \cup_{i=0}^{[\frac g2]}E_i^m$. Here $E_i^m$ is given similarly as \eqref{DMbdy}, which are essentially products of
lower dimensional moduli. We may assume $F_n(S)\subset E_i^m$ for some $0\leqslant i\leqslant [\frac g2]$. Finally, note that $\pi_n^{-1}(E_i^m)\subset \pd M_{g,n}^m$.
\end{cproof}

Now consider a  nonconstant holomorphic map $f:S\rightarrow \overline{M_g^m}$. By Theorem \ref{MT1}, $f(S)\subset \partial M^m_g=\cup_{i=0}^{[\frac{g}{2}]}E_{i}^m$ for $m$ sufficiently large.  
To simplify notation, we will drop the level `$m$' in $M^m_g, M^m_{g,n}, E^m_i$ etc, and suppose that $m$ is taken to be sufficiently large in the following.  
Assume $f(S)\subset E_i$ for some $0\leqslant i\leqslant [\frac g2]$. 
By the claim above, either 
$f(S)\subset \partial E_i$ or 
$f(S)\subset \pi_n^{-1}(x_0)$ for some $x_0 \in \pi_n(E_i-\partial E_i)$.

Suppose $f(S)$ does not lie on $\pd E_i$.  We consider each of the cases in (\ref{DMbdy}):

\ms
\noindent(i). $E_0 =\oM_{g-1,2}$.  In this case, $f(S)$ has to lie on a fiber of the forgetting map $\pi: M_{g-1,2}\rightarrow M_{g-1}$. 
For $x\in M_{g-1}$, the fiber $\pi^{-1}(x)\cong (R_x\times R_x-\Delta_x)/\mathbb{Z}_2$, where $R_x$ is the Riemann surface of genus $g-1$ represented by $x\in M_{g-1}$ and $\Delta_x$ is
the diagonal in $R_x\times R_x$.  Then $f(S)$ is
biholomorphic to a curve of genus $g_0$ in $(R_x\times R_x-\Delta_x)/\mathbb{Z}_2$, which lifts to a curve in $R_x\times R_x-\Delta_x$ since the quotient is unramified.  Note that for the
case that $f$ is an embedding, $g_0\geqslant g-1$ by Riemann-Hurwitz Formula.
 
\ms
\noindent
(ii). $E_i \cong \oM_{i,1}\times \oM_{g-i,1}, 1\leqslant i<\frac g2$.  In such case, the same argument as above shows that $f(S)$ has to lie on a fiber of the direct product forgetting map $\pi:M_{i,1}\times M_{g-i,1}\rightarrow M_{i}\times M_{g-i}$.  Let $x\in  M_{i}\times M_{g-i}$.  Then $\pi^{-1}(x)\cong R_x\times T_x$, where $R_x$
is a Riemann surface of genus $i$ and $T_x$ is a Riemann surface of genus $g-i$.  Again, if $f$ is an embedding, $g_0\geqslant \min(i,g-i)$ by Riemann-Hurwitz Formula.

\ms
\noindent
(iii). $E_{\frac g2} \cong(\oM_{\frac g2,1}\times \oM_{\frac g2,1})/\mathbb{Z}_2$, if $2|g$. In such case, the same argument as above shows that $f(S)$ has to lie on a fiber of the projection
$$\pi: M_{\frac g2,1}\times M_{\frac g2,1}\rightarrow (M_{\frac g2,1}\times M_{\frac g2,1})/\mathbb{Z}_2
\rightarrow M_{\frac g2}\times M_{\frac g2}.$$  
Let $x\in  M_{\frac g2}\times M_{\frac g2}$.  Then $\pi^{-1}(x)\cong R_x\times T_x$, where $R_x, T_x$
ares Riemann surface of genus $\frac g2$.  Again, if $f$ is an embedding, $g_0\geqslant \frac g2$ by Riemann-Hurwitz Formula.

\ms
Combining the three cases, we obtain the statement (II) and we are done. Therefore it remains to consider the case $f(S)\subset \partial E_i$. Note that the components of $\partial E_i$ is a direct product of at most two factors (possibility with mod $\mathbb{Z}_2)$, with factors of the from $\partial M_{g',n'}$ or $ M_{g',n'}$ for some $g'<g, n'<n$. By replacing $\pi_n$ in the claim above with the $\pi$ in each of the above three cases, we obtain a similar conclusion that $f(S)$ has to lie either on the boundary of a component of $\partial E_i$, or $f(S)$ lies in the fibre of $\pi$.
Hence by repeating the argument inductively, the image $f(S)$ has to lie on a stratum of the boundary.
\end{proof}

\section{Appendix: Sharp volume lower bounds}
The following function taken from \cite[Proposition 3.1.2]{HT2000} allows us to obtain a sharper volume estimate similar to that of \cite{HT2002}.
\begin{proposition}
Let $\varepsilon>0, 0<t_0<1$ be given real numbers and $N\in \mathbb{N}$ be a fixed positive integer. Then there exists a function $\psi_\varepsilon:(0,\infty)\rightarrow \mathbb{R}$ such that: \\
i) $\psi_\varepsilon\in \mathcal{C}^2((0,\infty))$ and $\psi_\varepsilon(t)= 0$ for $t\geq t_0$;\\
ii) $\psi_\varepsilon(t)$ is an increasing function for $t\in(0,\infty)$;\\
iii) $\psi_\varepsilon'(t)+t\psi_\varepsilon''(t)\geq - \frac{1}{(1-t)^2}$ for $t\in (0,\infty)$; and\\
iv) $\frac{t_0}{1-t_0}\geq 
\lim_{t\rightarrow 0} \frac{\psi_\varepsilon(t)}{\log t}\geq \frac{t_0}{1-t_0}-\varepsilon$.
\label{cutoff}
\end{proposition}
As in \cite{HT2002}, the volume lower bound for subvarieties in the case of Hermitian locally symmetric spaces obtained using Proposition \ref{cutoff} is sharp in the sense that such bound is realised by some totally geodesic subvarieties. In our case for obtaining sharper volume estimate, we may also apply Proposition \ref{cutoff} to construct another function $\phi$ as in the proof of Lemma \ref{lemvolest}: 
\begin{lemma}
Let $\phi(y)=\psi_\varepsilon(r(y))$. Then $\sqrt{-1}\ddbar \phi +\tomega\geqslant 0 $.
\end{lemma}
\begin{proof}
Write $\ell=\ell_C, r=r_C=\tanh \ell$ . We drop $i$ in below for convenience. By direction computation, $\barD r = (\sech^2\ell)\barD\ell=(1-r^2)\barD\ell$, so that
\begin{align*}
\ddbar r &= (1-r^2)\ddbar \ell+(-2r\partial r)\wedge \barD \ell \\
	&= \bigg(\frac{1-r^4}{r}-2r(1-r^2)\bigg) \partial\ell\wedge \barD\ell 
	\tag{$\because \ddbar\ell=\frac{1+r^2}{r}\partial\ell\wedge\barD\ell$.}\\
	&=\frac{(1-r^2)^2}{r}\partial\ell\wedge \barD\ell
	= \frac{1}{r}\partial r\wedge \barD r.
\end{align*}
Here $\ddbar\ell=\frac{1+r^2}{r}\partial\ell\wedge\barD\ell$ is interpreted as current  as in the proof of Lemma \ref{rcproperty}.
Therefore
\begin{align*}
\sqrt{-1}\ddbar \phi
&= \sqrt{-1}\psi'(r)\ddbar r + \sqrt{-1}\psi''(r)\partial r\wedge \barD r \\
&= \frac{1}{r}\bigg( \psi'(r) + r \psi''(r) \bigg) \sqrt{-1}\partial r\wedge \barD r\\
&\geqslant \bigg(-\frac{1}{r(1-r)^2}\bigg)\sqrt{-1}\partial r\wedge \barD r 
\tag{by Proposition \ref{cutoff} (iii)}\\
&\geqslant \bigg(-\frac{1}{(1-r)^2}\bigg)\sqrt{-1}\partial r\wedge \barD r  \\
&= -\tomega.
\end{align*}
\end{proof}

\section*{Declarations}
\textbf{Data availability statement:}
Not applicable to this article as no datasets were generated or analysed during
the current study.\\

\noindent
\textbf{Conflict of interest:} On behalf of all authors, the corresponding author states that there is no conflict of interest.

\end{document}